\documentclass{article}
\usepackage{graphicx}
\usepackage{amsmath,amssymb,amsfonts,latexsym}

\usepackage{natbib}
\usepackage{epic} 

\newtheorem{theorem}{Theorem}[section]
\newtheorem{lemma}[theorem]{Lemma}

\newcommand{\beq}{\begin{equation}}
\newcommand{\eeq}{\end{equation}}

\newcommand\Xtil{\widetilde{X}}
\newcommand\ybar{\bar{y}}
\newcommand\e{\varepsilon}

\DeclareMathOperator{\sign}{sign}

\newcommand\cF{\mathcal{F}}

\newcommand\cK{\mathcal{K}}

\newcommand\bE{\mathbf{E}}
\newcommand\bvar{\mathbf{Var}}
\newcommand\bP{\mathbf{P}}
\newcommand\bQ{\mathbf{Q}}

\newcommand\NN{\mathbb{N}}
\newcommand\RR{\mathbb{R}}
\newcommand\ZZ{\mathbb{Z}}

\newcommand\hgt{h}  
\newcommand\hgti{h}  
\newcommand\hgty{z}  
\newcommand\hgtyi{z}  
\newcommand\stsp{\mathcal{X}}  
\newcommand\hrap{\sigma} 
\newcommand\hham{z} 

\newcommand{\abs}[1]{\lvert#1\rvert}

\begin{document}

\title{Directed random growth models on the plane}
\author{Timo Sepp\"al\"ainen
\footnote{Research supported in part by NSF Grant DMS-0701091.}}

\maketitle


\begin{abstract} This is a brief survey of laws of
large numbers, fluctuation results and large deviation 
principles for asymmetric interacting
particle systems  that represent
moving interfaces on the plane.  We discuss the exclusion process,  
 the Hammersley process and the related last-passage growth
models. 
 \end{abstract}

\section{Introduction}
This article is a brief overview of recent results for a class of
stochastic processes that represent growth or motion of an 
interface in two-dimensional Euclidean space.   The models discussed have
in a sense rather orderly evolutions, and the word ``directed'' 
 is included in the title to evoke this feature.  

Let us begin with generalities about these stochastic processes.
The state at time $t\in[0,\infty)$ is
of the form  $\hgt(t)=(\hgti_i(t): i\in\ZZ)$ 
with the interpretation that the integer- or real-valued  
random variable $\hgti_i(t)$ represents the height of the 
interface over site $i$ of the substrate $\ZZ$.  
We call the state $\hgt=(\hgti_i)$ a {\sl height function} on $\ZZ$.  
The interface on the plane is then  represented by the graph 
$\{(i,\hgti_i): i\in\ZZ\}$. 
Each particular process has a state space that defines the set of admissible
height functions.  The state space will be defined by putting 
restrictions on the  increments (discrete derivatives)
 $\hgti_i-\hgti_{i-1}$ of the height 
functions.  

The random dynamics of the state are specified by the {\sl jump rates}
of the individual height variables $\hgti_i$.  
The rates are functions of the current state $\hgt$. 
For the sake of illustration, suppose  that only jumps 
 $\pm 1$ are permitted for each variable $\hgti_i$. 
Then the model is defined by giving two functions 
$p(\hgt)$ and $q(\hgt)$. If the current state is $\hgt$, then
the height value $\hgti_0$ at the origin jumps down to a new
value $\hgt_0'=\hgt_0-1$ with rate $p(\hgt)$, and  
 jumps up to a new
value $\hgt_0'=\hgt_0+1$ with rate $q(\hgt)$. 
In the spatially
homogeneous case  the rates for $\hgti_i$ are 
$p(\theta_i\hgt)$ and $q(\theta_i\hgt)$ where the spatial 
translations are defined by $(\theta_i\hgt)_j=\hgti_{i+j}$. 

The quantity   $p(\hgt)$  is the rate
for the jump  $\hgti_0\curvearrowright \hgti_0-1$ in an 
instantaneous sense: in the current state $h$,  the 
probability that this jump happens in the next infinitesimal 
time interval $(0, dt)$
  is $p(\hgt)dt+o(dt)$.  Rigorous constructions of the processes
utilize  Poisson processes or  ``Poisson clocks.'' 
A rate $\lambda$ Poisson process $N(t)$ is a simple
continuous time Markov chain: it starts at $N(0)=0$, runs
through the integers $0,1,2,3,\dotsc$ in increasing order, and waits 
for a rate $\lambda$ exponential random time between jumps. 
A  rate $\lambda$ exponential random time is defined by its
 density   $\varphi(t)=\lambda e^{-\lambda t}$ on $\RR_+$.  
The number of jumps $N(s+t)-N(s)$ in 
time interval $(s,s+t]$  has the mean $\lambda t$  Poisson  distribution
\[
P\{N(s+t)-N(s)=k\}= \frac{e^{-\lambda t}(\lambda t)^k}{k!}
\quad (k\ge 0). \]
If the overall rates are bounded, say by $p(\hgt)\le \lambda$, 
then a Poisson clock with a time-varying rate $p(\hgt(t))$
can be obtained from  $N(t)$  by 
randomly accepting a jump at time $t$ with probability 
$p(\hgt(t))/\lambda$. 

Later 
 we mention in passing  rigorous constructions of some processes.
In each case the outcome of the construction is that 
all the random variables $\{\hgti_i(t):i\in\ZZ,\, t\ge 0\}$ 
are defined as measurable functions
on an underlying probability space $(\Omega,\cF,\bP)$.  Since these processes
evolve through jumps, the appropriate path regularity is that
with probability 1, the path $t\mapsto \hgt(t)$ is right-continuous
with left limits (cadlag for short).  The use of Poisson clocks   
 makes the stochastic
process $\hgt(t)$ a {\sl Markov process}.  This means that 
if the present state $\hgt(t)$ is known, the future
evolution $(\hgt(s): s>t)$ is statistically independent
of the past  $(\hgt(s): 0\le s<t)$.  This is a consequence of the 
  ``forgetfulness
property'' of the exponential distribution.  For a complete
discussion of these  foundational  matters we must refer
to textbooks on probability theory and stochastic processes.

{\sl This article  covers only asymmetric systems.}  
Asymmetry in this context means that the height variables 
$\hgti_i(t)$ on average tend to move more in one direction
than the other.  For definiteness, we define the models so 
that the downward direction is the preferred one.
In fact, the great majority of the paper is concerned with
{\sl totally asymmetric systems} for which $q(\hgt)\equiv 0$,
in other words only downward jumps are permitted.
Symmetric systems behave quite differently from asymmetric
systems, hence restricting treatment to one or the other 
is natural.    

  Stochastic processes with
a large number of interacting components such as the height
process 
$\hgt(t)=(\hgt_i(t):i\in\ZZ)$ 
belong in an area of probability theory 
called {\sl interacting particle systems}.  
\citep{spit-70}  is one of 
 the seminal papers of this subject. 
  Here is a selection of 
books and lecture notes on the topic:
\citet{dema-pres}, \citet{durr-88}, \citet{kipn-land}, \citet{ligg85},
\citet{ligg99},  \citet{ligg04}, \citet{vara}.   
 \citet{krug-spoh} and \cite{spoh} are sources that combine 
 mathematics and the theoretical 
physics side. 

\medskip

Our treatment 
 is organized around  three basic 
questions  posed about stochastic models:
(i) laws of large numbers, (ii) fluctuations, and (iii) large deviations. 
 
\medskip

{\sl (i) Laws of large numbers} give deterministic
 limit shapes and evolutions
under appropriate space and time scaling. 
 A parameter
$n\nearrow\infty$ gives the ratio of macroscopic and microscopic
 scales.   A sequence of 
 processes $\hgt^n(t)$ indexed by $n$ is considered. 
Under appropriate hypotheses the height   process  satisfies this type of 
result:  for $(t,x)\in\RR_+\times\RR$  
\beq
n^{-1}\hgt^n_{[nx]}(nt)\to u(x,t) \quad\text{as $n\to\infty$,}
\label{lln1}\eeq
and the limit function $u$ satisfies a Hamilton-Jacobi equation
$
u_t+f(u_x)=0.
$

\medskip

{\sl (ii) Fluctuations.}
After  a law of large numbers 
the next question concerns  the 
random fluctuations around the large scale behavior. 
One seeks an exponent $\alpha$ that describes the magnitude
of these fluctuations, and a precise description of them
in the limit.  A typical statement would be
\beq
\frac{\hgti^n_{[nx]}(nt)- nu(x,t)}{n^\alpha}
\longrightarrow Z(t,x)
\label{fluct1}\eeq
where $Z(t,x)$ is a random variable whose distribution 
would be described  as part of
the result.  The convergence is of a weak type, 
where it is the probability distribution of the random variable 
on the left that converges. 

\medskip

{\sl (iii)  Large deviations.}  
 The vanishing probabilities of 
atypical behavior fall under this rubric. 
Often these  probabilities decay as
 $e^{-Cn^\beta}$ to leading order,
 with another exponent $\beta>0$.
The precise constant $C\in(0,\infty)$ is also of interest and comes in
the form of a {\sl rate function}. 
 When all the ingredients are in place the result is called  
a large deviation principle (LDP).  An LDP from the law of
large numbers \eqref{lln1} with rate function 
$I:\RR\to[0,\infty]$ 
 could take this form:
\beq
\lim_{\e\searrow 0}\lim_{n\to\infty}
n^{-\beta}\log \bP\{ \hgt^n_{[nx]}(nt)\in (z-\e,z+\e)\}
=-I(z)
\label{ldp1}\eeq
 valid for points $z$ in some range.  Positive
values $I(z)>0$ represent atypical behavior, while  
 limit
\eqref{lln1} would force $I(u(x,t))=0$.

\medskip

{\bf Example.} 
For classical examples of these statements let us
consider one-dimensional nearest-neighbor random walk.
Fix a parameter $0<p<1$. 
Let $\{X_k\}$ be independent, identically distributed  (IID)
$\pm 1$-valued  random
variables with common distribution 
$\bP\{X_k=1\}=p= 1- \bP\{X_k=-1\}$.   Define the 
random walk by $S_0=0$, $S_n=S_{n-1}+X_n$ for $n\ge 1$. 
Then the {\sl strong law of large numbers} gives the long term velocity:
\[
\lim_{n\to\infty} n^{-1}S_n= v \quad\text{where $v=\bE(X_1)=2p-1$.}
\]
The convergence in the limit above is almost sure (a.s.), that is,
almost everywhere (a.e.)  on the underlying probability space of 
the variables $\{X_k\}$. 
 
The order of nontrivial fluctuations around the limit is 
$n^{1/2}$ (``diffusive'') and in the limit these 
fluctuations are Gaussian.
That is the content of the {\sl central limit theorem}:
\[
\lim_{n\to\infty}
\bP\Bigl\{ \frac{S_n-nv}{\sigma n^{1/2}}\le s\Bigr\} 
=\Phi(s) \equiv \frac1{\sqrt{2\pi}}\int_{-\infty}^s 
e^{-z^2/2}\,dz.
\]
The parameter  $\sigma^2=\bE[(X_1-v)^2]$ is the variance. 

Random walk satisfies this LDP: 
\beq
\lim_{\e\searrow 0}\lim_{n\to\infty}
\frac1n\log \bP\{ \lvert S_n-nx\rvert\le n\e\} 
= -I(x) 
\label{rwldp1}\eeq
where the  rate function $I \colon\RR\to[0,\infty]$ is identically
$\infty$ outside $[-1,1]$ and  
\beq
I(x)= \frac{1-x}2\log \frac{1-x}{2(1-p)} 
+  \frac{1+x}2\log \frac{1+x}{2p} \quad\text{for $x\in[-1,1]$.}
\label{rwldp2}\eeq
$I(x)$ can be interpreted as an entropy.  Convex analysis 
plays a major role in large deviation theory.  Part of the 
general theory behind this simple case is that $I$ is
the convex dual of the logarithmic moment generating
function $\Lambda(\theta)=\log \bE(e^{\theta X_1})$. 

Results for random walk are covered in
graduate probability texts such as \citep{durr} and
\citep{kall}. 

\medskip

At the outset we delineated the class of models discussed.
Important models left out include 
{\sl diffusion limited aggregation} (DLA) and  
{\sl first-passage percolation}.  
Their interfaces
 are considerably more complicated than interfaces 
described by  height functions. 
But even for the models discussed our treatment is not
a complete  representation of the mathematical progress of
the past decade.  In particular, this article does not
delve into the recent work on Tracy-Widom fluctuations,
 Airy processes and determinantal point processes.  These
topics  are covered by many authors, and   we give a number of
references to the literature in Sections \ref{asepsec} 
and \ref{hammsec}. 
Overall, the best hope for this article is that it might inspire 
the reader to look further into the references. 

{\sl Recurrent notation.}  The set of nonnegative
integers is $\ZZ_+=\{0,1,2,\dotsc\}$, while $\NN=\{1,2,3,\dotsc\}$.
The integer part of a real $x$ is $[x]=\max\{n\in\ZZ: n\le x\}$. 
$a\vee b=\max\{a,b\}$ and $a\wedge b=\min\{a,b\}$.

{\sl Acknowledgement.}  The author thanks M.~Bal\'azs for 
valuable consultation during the preparation of this article. 

\section{Limit shape and evolution}
\label{llnsec}

\begin{figure}[ht]
\begin{center}
\begin{picture}(220,140)(5,20)
\put(40,20){\vector(1,0){140}}
\put(40,20){\vector(0,1){119}}
\multiput(40,40)(0,20){5}{\line(1,0){125}}
\multiput(63,20)(23,0){5}{\line(0,1){110}}
\put(44,26){$Y_{1,1}$}
\put(67,26){$Y_{2,1}$}
\put(90,26){$Y_{3,1}$}
\put(44,46){$Y_{1,2}$}
\put(67,46){$Y_{2,2}$}
\put(44,66){$Y_{1,3}$}
\put(184.7,17){$i$} \put(30,135){$j$}
\put(50,10){\small 1}
\put(73,10){\small 2}
\put(96,10){\small 3}
\put(119,10){\small 4}
\put(142,10){\small 5}
\put(29,27){\small 1} \put(29,47){\small 2} \put(29,67){\small 3}
\put(29,87){\small 4} \put(29,107){\small 5}
\end{picture}
\end{center}
\caption{Each point $(i,j)\in\NN^2$ has a weight $Y_{i,j}$ attached
to it.} \label{fig1}
\end{figure}
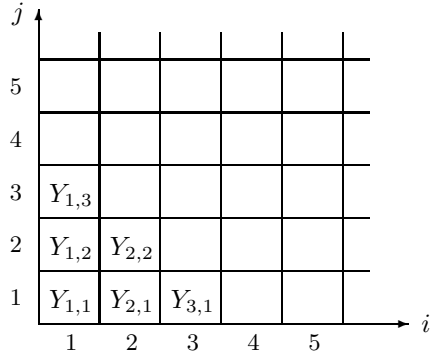

We begin with the much studied  {\sl corner
growth model} and a description that is not directly 
in terms of height variables. 
Attach 
nonnegative weights  $\{Y_{i,j}\}$ to the points $(i,j)$ 
of the positive quadrant $\NN^2$ of $\ZZ^2$, as in Figure \ref{fig1}. 
$Y_{i,j}$ represents the time it takes to occupy point $(i,j)$ 
{\sl after} the points to its left and below have been 
occupied.    Assume that everything outside the positive
quadrant is occupied at the outset so the process can start.   
Once occupied, a point remains occupied.  Thus this is  a
totally asymmetric growth model, for the growing cluster never 
loses points, only adds them.

\begin{figure}[ht]
\begin{center}
\begin{picture}(220,140)(5,20)
\put(40,20){\vector(1,0){140}}
\put(40,20){\vector(0,1){119}}
\multiput(40,40)(0,20){5}{\line(1,0){125}}
\multiput(63,20)(23,0){5}{\line(0,1){110}}
\multiput(49,27)(23,0){3}{\large$\bullet$}
\multiput(95,47)(0,20){2}{\large$\bullet$}
\multiput(118,67)(23,0){2}{\large$\bullet$}
\put(141,87){\large$\bullet$} 
\multiput(56.5,30)(23,0){2}{\vector(1,0){14}}
\multiput(97.8,33.8)(0,20){2}{\vector(0,1){12}}
\multiput(102.5,70)(23,0){2}{\vector(1,0){14}}
\put(144,73.8){\vector(0,1){12.2}}
\put(184.7,17){$i$} \put(30,135){$j$}
\put(50,10){\small 1} \put(73,10){\small 2}
\put(96,10){\small 3} \put(119,10){\small 4}
\put(142,10){\small 5}
\put(29,27){\small 1} \put(29,47){\small 2} \put(29,67){\small 3}
\put(29,87){\small 4} \put(29,107){\small 5}
\end{picture}
\end{center}
\caption{An admissible path from $(1,1)$ to $(5,4)$}  \label{fig2}
\end{figure}
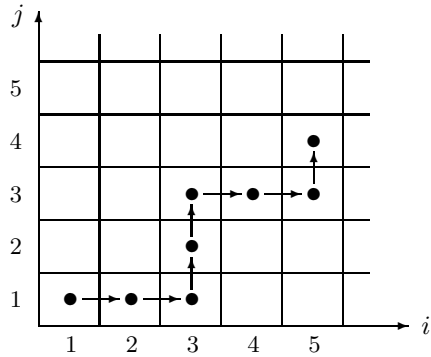

Let $G(k,\ell)$ denote the time when point $(k,\ell)$ becomes
occupied.  The above explanation is summarized 
by these rules: 
$G(k,\ell)=0$ for $(k,\ell)\notin\NN^2$, and  
\beq G(k,\ell)= G(k-1,\ell)\vee G(k,\ell-1)+ Y_{k,\ell}
\quad\text{ for $(k,\ell)\in\NN^2$.}
\label{Grec}\eeq
  The last equality can be iterated until
the corner $(1,1)$ is reached, resulting in this last-passage
formula for $G$:
\beq
G(k,\ell)=\max_{\pi\in\Pi_{k,\ell}} \sum_{(i,j)\in\pi} Y_{i,j}
\label{Geq1}\eeq 
where $\Pi_{k,\ell}$ is the collection of nearest-neighbor up-right
paths $\pi$ from $(1,1)$ to $(k,\ell)$.   Figure \ref{fig2}
represents one such path for $(k,\ell)=(5,4)$. 

This model and others of its kind are called 
 {\sl directed last-passage percolation models}. 
``Directed'' refers to the restrictions on admissible paths,
and ``last-passage'' to the
feature that the occupation time $G(k,\ell)$ is determined 
by the slowest path to $(k,\ell)$.   (By contrast, in first-passage
percolation occupation times are determined by quickest paths.)

Our first goal is to argue the existence of a 
limit for $n^{-1}G([nx],[ny])$ as $n\to\infty$. 
Assume now that 
the weights $\{Y_{i,j}\}$ are IID nonnegative 
random variables. 

The idea is to exploit sub(super)additivity.  Generalize 
the definition of $G(k,\ell)$ to 
\[
G((k,\ell),(m,n))=\max_{\pi\in\Pi_{(k,\ell),(m,n)}} 
\sum_{(i,j)\in\pi} Y_{i,j},
\]
where $\Pi_{(k,\ell),(m,n)}$ is the collection of nearest-neighbor up-right
paths $\pi$ from $(k+1,\ell+1)$ to $(m,n)$.  
The definitions lead to
the superadditivity
\beq
G(k,\ell) +G((k,\ell),(m,n))\le G(m,n).
\label{supadd}
\eeq
Kingman's 
subadditive ergodic theorem \citep[Chapter 6]{durr}
 and some estimation 
 implies the existence of a 
deterministic limit 
\beq \gamma(x,y)=\lim_{n\to\infty}n^{-1} G([nx],[ny]) 
\quad\text{for $(x,y)\in\RR_+^2$.}
\label{Glim}\eeq

Moment assumptions under which the limit function
$\gamma$ is continuous up to the boundary were investigated
by \citet{mart-04}.   We turn to the problem of computing
$\gamma$ explicitly, and for this  we need very specialized 
assumptions.  Essentially only one distribution can be currently 
handled: the exponential, and its discrete counterpart, 
the geometric.   Take the $\{Y_{i,j}\}$ to be  IID rate 1 exponential 
random variables. In other words their common density is
$e^{-y}$. 

The difficulty with finding the explicit  limit has to do
with the superadditivity.  The limit in 
Birkhoff's ergodic theorem is simply the 
expectation of the function averaged over shifts:
$n^{-1}\sum_{k=1}^n f\circ\theta_k\to \bE f$
\citep[Chapter 6]{durr}.  But the 
subadditive ergodic theorem gives only an asymptotic
expression for the limit.   
We need a new ingredient.  We shall embed the last-passage
model into the totally asymmetric simple exclusion process
(TASEP). This has 
explicitly identifiable invariant distributions (``steady states'')
 with which we can do explicit calculations.

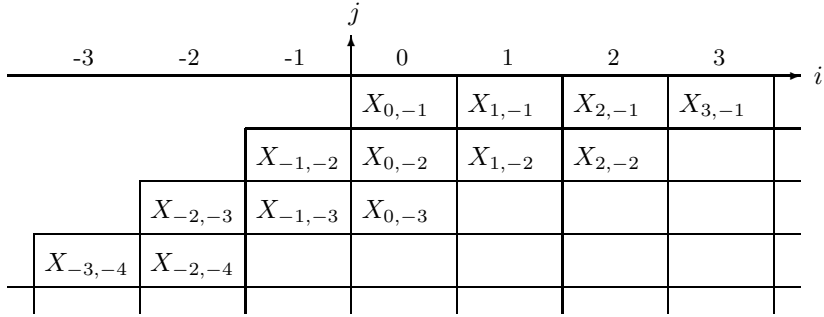
\begin{figure}[ht]
\begin{center}
\begin{picture}(300,120)(-130,-80)
\put(0,0){\vector(1,0){170}} \put(175,-3){$i$} \put(-1,21){$j$}
\put(0,0){\line(-1,0){130}}
\put(0,0){\vector(0,1){15}}
\put(0,0){\line(0,-1){90}}
\multiput(40,0)(40,0){4}{\line(0,-1){90}} 
\put(-40,-20){\line(0,-1){70}\line(1,0){210}} 
\put(-80,-40){\line(0,-1){50}\line(1,0){250}} 
\put(-120,-60){\line(0,-1){30}\line(1,0){290}} 
\put(-130,-80){\line(1,0){300}} 
\put(17,4){\small 0}\put(57,4){\small 1}\put(97,4){\small 2}
\put(137,4){\small 3} 
\put(-25,4){\small -1}\put(-65,4){\small -2}
\put(-105,4){\small -3}
\put(4,-14){$X_{0,-1}$}\put(44,-14){$X_{1,-1}$}\put(84,-14){$X_{2,-1}$}
\put(124,-14){$X_{3,-1}$}
\put(-36,-34){$X_{-1,-2}$}\put(4,-34){$X_{0,-2}$}\put(44,-34){$X_{1,-2}$}
\put(84,-34){$X_{2,-2}$}
\put(-76,-54){$X_{-2,-3}$}\put(-36,-54){$X_{-1,-3}$}\put(4,-54){$X_{0,-3}$}
\put(-116,-74){$X_{-3,-4}$}\put(-76,-74){$X_{-2,-4}$}
\end{picture}
\end{center}
\caption{Last-passage model for TASEP.  The horizontal
$i$-axis and the vertical $j$-axis are labeled, and points
on the $i$-axis from $-3$ to $3$ are marked. 
Weights $X_{i,j}$ are 
attached to points $(i,j)$ such that $i\in\ZZ$,
$j\in-\NN$ and  $j<0\wedge i$.}\label{fig3} 
\end{figure}%

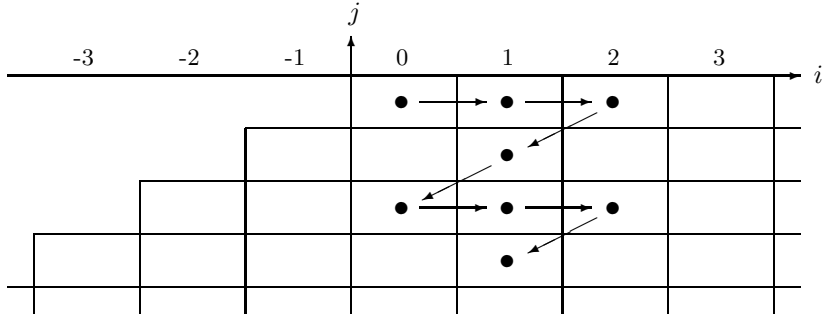
\begin{figure}[ht]
\begin{center}
\begin{picture}(300,120)(-130,-80)
\put(0,0){\vector(1,0){170}} \put(175,-3){$i$} \put(-1,21){$j$}
\put(0,0){\line(-1,0){130}}
\put(0,0){\vector(0,1){15}}
\put(0,0){\line(0,-1){90}}
\multiput(40,0)(40,0){4}{\line(0,-1){90}} 
\put(-40,-20){\line(0,-1){70}\line(1,0){210}} 
\put(-80,-40){\line(0,-1){50}\line(1,0){250}} 
\put(-120,-60){\line(0,-1){30}\line(1,0){290}} 
\put(-130,-80){\line(1,0){300}} 
\put(17,4){\small 0}\put(57,4){\small 1}\put(97,4){\small 2}
\put(137,4){\small 3} 
\put(-25,4){\small -1}\put(-65,4){\small -2}
\put(-105,4){\small -3}
\multiput(16,-13)(40,0){3}{\large $\bullet$} \put(56,-33){\large $\bullet$}
\multiput(16,-53)(40,0){3}{\large $\bullet$}  \put(56,-73){\large $\bullet$}
\multiput(26,-10)(40,0){2}{\bf\vector(1,0){25}}
\multiput(93,-14)(-40,-20){2}{\bf\vector(-2,-1){26}}
\multiput(26,-50)(40,0){2}{\bf\vector(1,0){25}}
\put(93,-54){\bf\vector(-2,-1){26}}
\end{picture}
\end{center}
\caption{The image of the path in Figure \ref{fig2}. Now it
goes from $(0,-1)$ to $(1,-4)$.}\label{fig4} 
\end{figure}

Originally TASEP was introduced as  a particle model.
We wish to link TASEP with the last-passage model in a
way that preserves the  original formulation of TASEP,
while mapping   particle occupation variables into height
increments  and  particle current into  column growth.
To achieve this  we
transform the coordinates $(i,j)$ of Figure \ref{fig1} 
via the bijection $(i,j)\mapsto (i-j,-j)$.
The result is the last-passage model of Figure \ref{fig3}. 
Weights are  relabeled as $X_{i,j}=Y_{i-j,-j}$.  
The transformation of admissible paths is illustrated 
by Figure \ref{fig4}.  Let the new last-passage times 
 be denoted by 
$H(k,\ell)$. 
 For $\ell<0\wedge k$ the maximizing-path formulation uses now
paths of the kind represented in Figure \ref{fig4}: 
\beq
H(k,\ell)=\max_{\sigma\in\Sigma_{k,\ell}} \sum_{(i,j)\in\sigma} X_{i,j}
\label{Heq1}\eeq 
where $\Sigma_{k,\ell}$ is the collection of 
paths $\sigma$ from $(0,-1)$ to $(k,\ell)$ that take
steps of two types: $(1,0)$ and $(-1,-1)$. 
  The connection with
the previous last-passage process is $H(k,\ell)=G(k-\ell,-\ell)$.
The process $\{H(k,\ell)\}$ is also defined by the 
 recursion 
\beq
H(k,\ell)= H(k-1,\ell)\vee H(k+1,\ell+1) + X_{k,\ell}
\quad(\ell<0\wedge k)
\label{Hrec}\eeq
together with the boundary values $H(k,\ell)=0$ for 
$\ell\ge 0\wedge k$.  It satisfies the 
limit
\beq
 \lambda(x,y)=\lim_{n\to\infty} n^{-1} H([nx],[ny]) = \gamma(x-y,-y)
\quad\text{for $y<0\wedge x$.}\label{Hlim}\eeq

To establish the TASEP connection
 we first define a height process $w(t)=(w_i(t):i\in\ZZ)$ that 
will turn out to be an alternative description of 
 the last-passage process $\{H(i,j)\}$.
Initially at time $t=0$ the height is given by 
\beq w_i(0)= \begin{cases}  i, &i\le -1\\
                         0, &i\ge 0. \end{cases}
\label{wedge0}\eeq
This is the boundary of the region $\{j<0\wedge i\}$
 filled with $X_{i,j}$'s 
in Fig.\ \ref{fig3}.   This initial shape is a wedge, hence the
symbol $w$.  

\begin{figure}[ht]
\begin{center}
\begin{picture}(300,120)(-130,-80)
\put(0,0){\line(1,0){170}} \put(175,-3){$i$} \put(-1,21){$j$}
\put(0,0){\line(-1,0){130}}
\put(0,0){\vector(0,1){15}}
\put(0,0){\line(0,-1){90}}
\multiput(40,0)(40,0){4}{\line(0,-1){90}} 
\put(-40,-20){\line(0,-1){70}\line(1,0){210}} 
\put(-80,-40){\line(0,-1){50}\line(1,0){250}} 
\put(-120,-60){\line(0,-1){30}\line(1,0){290}} 
\put(-130,-80){\line(1,0){300}} 
\put(17,4){\small 0}\put(57,4){\small 1}\put(97,4){\small 2}
\put(137,4){\small 3} 
\put(-25,4){\small -1}\put(-65,4){\small -2}
\put(-105,4){\small -3}
{\linethickness{3pt}
\put(-120,-58.6){\line(0,-1){21}} \put(-120,-60){\line(1,0){121}} 
\put(-130,-80){\line(1,0){11.3}}
\put(0,-38.7){\line(0,-1){22.7}} \put(0,-40){\line(1,0){81.3}} 
\put(80,-19){\line(0,-1){21}} \put(78.7,-20){\line(1,0){42.8}}
\put(120,0.3) {\line(0,-1){21}} \put(118.5,0){\line(1,0){51.7}}
}
\put(16,-53){$\times$}  \put(136,-13){$\times$} 
\put(-104,-73){$\times$} 
\end{picture}
\end{center}
\caption{A possible height function $w(t)$ (thickset
graph) with column values $w_{-1}(t)=-3$, 
$w_0(t)=-2$, $w_1(t)=-2$, etc.   8 jumps have
taken place during time $(0,t]$.   The columns grow
downward.  $\times$'s mark the
allowable jumps from this state.} \label{fig5}
\end{figure}
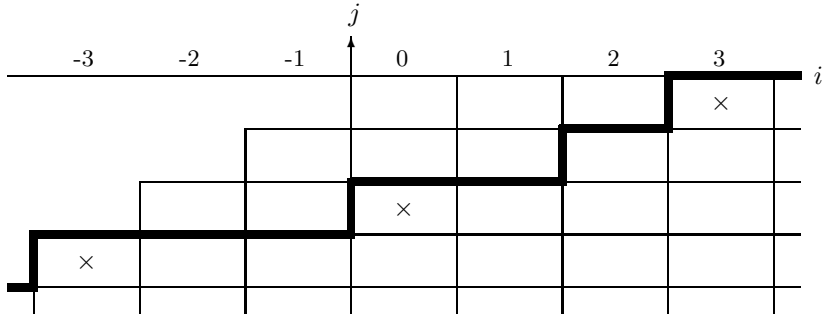

Give each column $i$ an independent rate 1 Poisson clock $N_i$. 
Variable $w_i$ jumps downward according to this rule: 
if $t$ is a jump time for Poisson process $N_i$, then 
\beq
\text{$w_i(t)=w_i(t-)-1$ provided }
\begin{cases} w_{i-1}(t-)=w_{i}(t-)-1 \text{  and }\\
  w_{i+1}(t-)=w_i(t-).  \end{cases} 
 \label{taseprule1}\eeq
Equivalently, each column variable $w_i$ jumps down independently at 
rate 1, as long as the state $w(t)$ remains in the state space
\beq
\stsp_1=\{\hgt\in\ZZ^\ZZ:  \text{$\hgti_{i}-\hgti_{i-1}\in\{0,1\}$  
for all $i\in\ZZ$} \}.
\label{tasepX}\eeq 
(See Fig.\ \ref{fig5} for an example.)  The {\sl interaction} 
between the variables  
is encoded in rule \eqref{taseprule1}. It forces the time-evolution
of each variable $w_i$ to depend on the evolution of its neighbors. 
By contrast, without the interaction the variables $w_i$ would
simply march along as Poisson processes  independently of each other.

This construction defines the height process
$w(t)=(w_i(t):i\in\ZZ)$ for all times $t\in[0,\infty)$
 in terms of the family of Poisson clocks $\{N_i\}$. 
Given this process  $w(t)$ define the stopping times 
\beq
T(i,j)=\inf\{t\ge 0:  w_i(t)\le j\} 
\label{Tij}\eeq
that mark the time when column $i$ first reaches level $j$.
{\sl Stopping time} is a technical term for a random time
whose arrival can be verified  without looking into the future.

Initial condition \eqref{wedge0}  implies $T(i,j)=0$ for
$j\ge i\wedge 0$. 
Rule \eqref{taseprule1} tells us that $T(i-1,j)\vee T(i+1,j+1)$ is
the stopping time at which the system is ready  for
$w_i$ to jump from level $j+1$  to $j$.  (Note that $w_i$ must have
reached level $j+1$ already earlier because if the rules are
followed, $T(i+1,j+1)$ comes after $T(i,j+1)$.) By the 
forgetfulness property of the exponential distribution, after the
stopping time   $T(i-1,j)\vee T(i+1,j+1)$
 it takes another independent rate 1 exponential 
time $\Xtil_{i,j}$ until $w_i$ jumps from $j+1$ to $j$.
Consequently the process $\{T(i,j)\}$ 
satisfies the recursion 
\beq
T(i,j)= T(i-1,j)\vee T(i+1,j+1) + \Xtil_{i,j}.  
\label{wedgerec}\eeq
This is of the same form as the recursion \eqref{Hrec} 
satisfied by  $\{H(i,j)\}$. 
   From this one can prove that indeed the
processes  $\{T(i,j)\}$ and  $\{H(i,j)\}$ are equal in distribution.
Therefore \eqref{Hlim} gives  also $n^{-1}T([nx],[ny])\to \lambda(x,y)$.
This is the precise meaning of the earlier claim that the height
process $w(t)$ gives an alternative description of the last-passage
process $\{H(i,j)\}$.   

Subadditivity and some estimation justifies the existence of 
a concave function $g$ on $\RR$ such that 
\beq
\lim_{t\to\infty} t^{-1}w_{[xt]}(t)=g(x)\quad\text{a.s.\ for $x\in\RR$.}
\label{wedgelim}\eeq
Since rates are $1$, $g$ records only the initial height outside 
the interval $[-1,1]$, and so  
\[
\text{$g(x)=0\wedge x$ for $\abs{x}>1$. }\] 
Since the interface is a level curve of passage times, 
$\lambda(x,g(x))=1$ for $-1\le x\le 1$.  The last-passage limits
are homogeneous in the sense that 
  $\lambda(cx,cy)=c\lambda(x,y)$ for $c>0$.
Consequently $\lambda$ and then $\gamma$ can be obtained from $g$.

\medskip

To summarize, 
thus far we have converted the original task of computing 
$\gamma(x,y)$ of \eqref{Glim} to finding the function $g$ 
of \eqref{wedgelim} on the interval $[-1,1]$. 
Now consider the general height process $\hgt(t)$ with state
space $\stsp_1$ from \eqref{tasepX}
 and dynamics defined as for $w(t)$ above:
height variable $h_i$ 
jumps one step down at every jump epoch of the Poisson
clock $N_i$, provided
this jump does not take the height function out of $\stsp_1$.  
A jump attempt that would violate the state space restriction
is simply ignored. 

Certain technical issues may trouble the reader. 
An infinite family of Poisson clocks has infinitely 
many jumps in any nonempty time interval $(0,\e)$. So  
there is no first jump attempt in the system and it is not 
obvious that the local rule leads to a well-defined global 
evolution:
to determine the evolution of $\hgti_i$ on $[0,t]$
 we need to look at the evolution of
its neighbors $\hgti_{i\pm 1}$, and then their neighbors, 
ad infinitum.  
However, given any $T<\infty$, almost surely there are   
 indices $i_k\searrow-\infty$ and $i_k'\nearrow\infty$ such that 
$N_{i_k}$ and $N_{i_k'}$ have no  jumps during $(0,T]$.  
Consequently the system decomposes into (random) finite pieces
that do not communicate before time $T$. The evolution 
can be determined separately in each finite segment which do 
experience only finitely many jumps up to time $T$ (again almost surely). 
Another technical point is that the clocks $\{N_i\}$ 
have no simultaneous jumps (almost surely)
 so one never needs to consider more
than one jump at a time. 

Given that the height process $\hgt(t)$ has been constructed, next
define the increment process $\eta(t)=(\eta_i(t):i\in\ZZ)$ by 
\beq
\eta_i(t)=\hgti_i(t)-\hgti_{i-1}(t). \label{defeta1}\eeq
Process $\eta(t)$ has compact state space $\{0,1\}^\ZZ$ and 
its dynamics inherited from $\hgt(t)$
 can be succinctly stated as follows: each {\tt 10}
pair becomes a {\tt 01} pair at rate 1, independently of
the rest of the system. To see this connection, 
observe that if Poisson clock $N_i$ jumps at time $t$,
the height process  undergoes the  transformation 
$\hgti_i(t)=\hgti_i(t-)-1$ 
only if $(\eta_i(t-)=1,\eta_{i+1}(t-)=0)$, and then after the 
jump the situation is $(\eta_i(t)=0,\eta_{i+1}(t)=1)$.
This is a direct translation of the condition that jumps
are executed only if the state $\hgt$ remains in the state space
$\stsp_1$.  

  It is natural to interpret 
the {\tt 1}'s as particles and the {\tt 0}'s as holes, or vacant sites. 
The process $\eta(t)$ is the {\sl totally asymmetric simple 
exclusion process} (TASEP). 
In this model the only interaction between the particles 
is the exclusion rule  that stipulates that
particles are not allowed to jump onto occupied sites.
This property is enforced  by the evolution
because the definitions made above ensure
that a jump in Poisson clock $N_i$ sends a particle
from site $i$ to site $i+1$ only if site $i+1$ is vacant.
Total asymmetry refers to the property that particles jump
only to the right, never left. 
The definitions also entail this connection between 
the heights and the particles:  
\beq
\text{$\hgti_i(0)-\hgti_i(t)$ $=$ cumulative particle current
across the edge $(i,i+1)$.}
\label{h-curr1}\eeq 

We need to discuss  two more properties of these processes,
(i) stationary behavior and (ii) the envelope property. 
Then we are ready to compute the function $g$ of 
\eqref{wedgelim}.

{\sl (i) Stationary behavior.}
For $\rho\in[0,1]$,  the {\sl Bernoulli} probability measure $\nu_\rho$ 
on $\{0,1\}^\ZZ$ is defined by the requirement that  
\beq
\nu_\rho\{\eta: \text{$\eta_i=1$ for $i\in I$, 
$\eta_j=0$ for $j\in J$} \} =\rho^{\abs{I}}(1-\rho)^{\abs{J}}
\label{defnurho}\eeq
for any disjoint $I,J\subseteq\ZZ$ with cardinalities 
$\abs{I}$ and $\abs{J}$.  Measure  $\nu_\rho$ corresponds
to putting a particle at each site independently
 with probability $\rho$. 

It is known that the measures $\{\nu_\rho\}_{\rho\in[0,1]}$
are invariant for the process $\eta(t)$, and in fact they are
the extremal members of the compact, convex set of 
invariant probability measures that are also invariant under spatial shifts. 
Invariance means that if the process $\eta(t)$ is started
with a random $\nu_\rho$-distributed initial state $\eta(0)$,
then at each time $t\ge 0$ the state $\eta(t)$ is 
$\nu_\rho$-distributed, and furthermore, the probability
distribution of the entire process
$\eta(\cdot)=(\eta(t):t\ge 0)$ is invariant under time shifts.  

If we know the current state $h(t)$, then the 
probability that  $h_0$ jumps down in a short time interval
$(t,t+\e)$  
is   $\e\eta_0(t)(1-\eta_1(t))+O(\e^2)$. 
This follows because a jump can happen only when a {\tt 10} pair
is present, and   from  properties
of Poisson processes.  
Estimation of this kind proves that
\beq
\hgti_0(t)-\hgti_0(0)=-\int_0^t \eta_0(t)(1-\eta_1(t))\,ds
+M(t)
\label{hmg}\eeq
where $M(t)$ is a mean-zero martingale. This
identity is a  stochastic ``fundamental theorem 
of calculus'' of sorts.     Since things are random the difference
between $\hgti_0(t)-\hgti_0(0)$ and the 
integral of the infinitesimal rate  cannot be 
 identically zero. Instead it is a {\sl martingale}.  This is a process
whose increments have mean zero in a very strong sense, namely even 
 when conditioned on the entire
past. 

Let us average over \eqref{hmg} in the stationary situation.
  Let $\bE_{\nu_\rho}$ denote expectation of functions of the
stationary  process $\eta(\cdot)$ whose state $\eta(t)$ 
is $\nu_\rho$-distributed at each time $t$.  Normalize the height 
process $\hgt(\cdot)$ at time zero
 so that $\hgti_0(0)=0$.  Then $\hgt(\cdot)$ is
entirely determined by $\eta(\cdot)$.
Since $\eta_0(s)$ and $\eta_1(s)$ are independent
at any fixed time $s$,   we get 
\beq E_{\nu_\rho}[\hgti_0(t)]= -tf(\rho) \label{Ehgt1}\eeq
 where the {\sl particle flux}
is defined by 
\beq f(\rho)=\rho(1-\rho). \label{def-f}\eeq

{\sl (ii) Envelope property.}    Even though the flux $f$ is
nonlinear and therefore, as we see later, TASEP is governed by
a nonlinear PDE,  the height process has a valuable additivity 
property.  Suppose a given initial height function
$\hgt(0)\in\stsp_1$  is the envelope of  a countable collection 
$\{\hgty^{(k)}(0)\}_{k\in\cK}$ of height functions
in the sense that 
\beq \hgti_i(0)=\sup_{k\in\cK} \hgtyi^{(k)}_i(0) 
\quad\text{for each site $i\in\ZZ$.}\label{env1}\eeq
  Take a single
collection $\{N_i\}$ of Poisson clocks, and let all processes
$\hgt(t)$, $\hgty^{(k)}(t)$ evolve from their initial height functions
 by following the same clocks
$\{N_i\}$. This kind of simultaneous construction of many random
objects for the purpose of comparison is called 
a {\sl coupling.} By induction on jumps one can prove that
this coupling preserves the envelope property
 for all time:

\begin{lemma} 
$\displaystyle{ \hgti_i(t)=\sup_{k\in\cK} \hgtyi^{(k)}_i(t)}$ 
for all $i\in\ZZ$ and all $t\ge 0$, almost surely. 
\label{varlm}\end{lemma}

We take as auxiliary processes 
 $\hgty^{(k)}(t)$ suitable translations
of the basic wedge process defined in 
\eqref{wedge0}--\eqref{taseprule1}.  For $k\in\ZZ$ set 
\beq w^{(k)}_i(0)  = \begin{cases}  i-k, &i< k\\
                         0, &i\ge k \end{cases}
\quad\text{and}\quad  \hgtyi^{(k)}_i(t)=\hgti_k(0)+w^{(k)}_i(t).
\label{wedgek}\eeq
The apex of the wedge $\hgty^{(k)}(0)$  is 
at the  point $(k, \hgti_k(0))$, and then the 
 definition of the wedge ensures that $\hgt(0)\ge\hgty^{(k)}(0)$. 
Hypothesis \eqref{env1} holds and
  Lemma \ref{varlm} gives this variational equality:
\beq
\hgti_i(t)=\sup_{k\in\ZZ}\bigl\{  \hgti_k(0)+w^{(k)}_i(t)\bigr\}.
\label{var1}\eeq

Now we extract two results from the assembled ingredients:
  the function $g$ of \eqref{wedgelim}, and a general ``hydrodynamic
limit'' that describes the large scale evolution of the process.

First specialize \eqref{var1} to the stationary situation where
$\hgti_0(0)=0$ and the increments  are $\nu_\rho$-distributed,
and write  \eqref{var1} in the form
\beq
t^{-1}\hgti_0(t)=\sup_{y\in\RR}\bigl\{ t^{-1} \hgti_{[ty]}(0)
+t^{-1}w^{([ty])}_0(t)\bigr\}.
\label{var2}\eeq
Let $t\to\infty$.  Inside the braces on the right 
$ t^{-1} \hgti_{[ty]}(0)\to \rho y$ a.s.\ by the law of large numbers.
 $t^{-1}w^{([ty])}_0(t)\to g(-y)$ by a translation of the limit 
\eqref{wedgelim}.  
With some work take the limit
outside the supremum.  Then we know $t^{-1}\hgti_0(t)$ converges.
By supplying some moment bounds we can take expectations over the limits,
and with  \eqref{Ehgt1}  arrive at 
\beq
-f(\rho)=\sup_{y\in\RR}\{ \rho y + g(-y)\}.
\label{dual1}\eeq
This is a convex duality
\citep{rock-ca}. From the explicit invariant distributions \eqref{defnurho}
 we obtained $f$ in \eqref{def-f}, and then  we can
solve \eqref{dual1}  for $g$. (Without the invariant
distributions we can carry out part of this reasoning but we cannot find
$f$ and $g$ explicitly.)   Let us record the results.

\begin{theorem}
For the 
limit {\rm\eqref{wedgelim}}
 $g(x)=-\tfrac14(1-x)^2$ for $-1\le x\le 1$. 
For the limit {\rm\eqref{Glim}}
 $\gamma(x,y)=(\sqrt{x}+\sqrt{y}\,)^2$ for $x,y\ge 0$.  
\label{gammathm} \end{theorem}

We turn to the hydrodynamic limit. 
Assume given a function $u_0$ on $\RR$
and a sequence of random initial height functions $\hgt^n(0)\in\stsp_1$
($n\in\NN$)  such that 
\beq
n^{-1}\hgti^n_{[nx]}(0)\to u_0(x) 
\quad\text{a.s.\ as $n\to\infty$ for each $x\in\RR$.}\label{llninit}\eeq 
For this to be possible $u_0$ has to be Lipschitz with
$0\le u_0'(x)\le 1$ Lebesgue-a.e.  

\begin{theorem}  For $x\in\RR$ and $t>0$ we have the limit
\beq
n^{-1}\hgti^n_{[nx]}(nt)\to u(t,x) 
\quad\text{a.s.\ as $n\to\infty$}\label{llnt>0}\eeq
where 
\beq
u(t,x)=\sup_{y\in\RR}\Bigl\{ u_0(y) +tg\Bigl(\frac{x-y}t\Bigr)\Bigr\}.
\label{hopflax}\eeq 
\label{hydrothm}
\end{theorem} 

Equation \eqref{hopflax} is a Hopf-Lax formula 
\citep{evan} and it says that  $u$ is the entropy solution of
the Hamilton-Jacobi equation 
\beq
u_t+f(u_x)=0,\quad u\vert_{t=0}=u_0. \label{hj1}\eeq
In other words this equation governs the macroscopic evolution of
the height process.   Theorem \ref{hydrothm} is proved by 
showing that, as $n\to\infty$,   variational formula \eqref{var1} for  
$n^{-1}\hgt^n_{[nx]}(nt)$ 
turns into \eqref{hopflax}.  Details can be found in 
\citep{sepp99K}.  

\medskip

{\bf Further remarks.} The function $g$ in Theorem \ref{gammathm}
 was first calculated by \citet{rost} in one of the 
seminal papers of hydrodynamic limits, but without the last-passage 
representation and with a different approach than
the one presented here.  

Let us discuss various avenues of generalization.
 We encounter immediately difficult open problems.   

\medskip

{\sl (i) Generalizations that retain the envelope property.}   
The argument sketched above that combines the envelope property
with the duality of the flux and the wedge shape
 to derive hydrodynamic limits
 was introduced in \citep{sepp98ebp, sepp98mprf,
sepp99K}.  An earlier instance of the variational connection
 appeared in 
 \citet{aldo-diac95} for Hammersley's process. 
This work itself was based on the classic paper \citep{hamm}; 
see Section \ref{hammsec} below.  Also in queueing
literature similar variational expressions arise 
 \citep{kell-szcz}.  

To define the {\sl $K$-exclusion process} we
 replace the state space $\stsp_1$ of \eqref{tasepX} 
with 
\beq
\stsp_K=\{\hgt\in\ZZ^\ZZ:  \text{$0\le \hgti_{i}-\hgti_{i-1}\le K$  
for all $i\in\ZZ$} \}
\label{KtasepX}\eeq 
for some $2\le K<\infty$.  Otherwise keep the model the same:
rate 1 Poisson clocks $\{N_i\}$ govern the jumps of height
variables $\hgti_i$, and jumps that take the state $\hgt$
outside the space $\stsp_K$ are prohibited.  
The increment process is now called totally asymmetric 
 $K$-exclusion (some authors use ``generalized exclusion'').   
The variational coupling (Lemma \ref{varlm}) works as before.
But  invariant distributions are unknown, and there is 
 even no  proof of existence of an extremal invariant distribution
for each density value $\rho\in[0,K]$.   No alternative way to 
compute $f$ and $g$  has been found.  Theorem \ref{hydrothm}
 is valid, but the 
most that can be said about $f$ and $g$ is that they exist as
concave functions.

Interestingly, the situation becomes again explicitly
 analyzable for $K=\infty$ where  the only constraint on $\hgt$ is  
$\hgti_i\le \hgti_{i+1}$.  The increment process is a special case
of a  {\sl zero range process}. Its state space
is $(\ZZ_+)^\ZZ$ and  IID geometric distributions are invariant
\citep{ligg-73, andj-82}.
As a final step of generalization, away from monotone height
functions, let us mention {\sl bricklayer
processes} \citep{bala-03, bala-etal-07} whose
 increments $\eta_i=\hgti_i-\hgti_{i-1}$
can be positive or negative.  

The variational coupling of Lemma \ref{varlm} works equally well
for certain multidimensional height processes
 $\hgt(t)=(\hgti_i(t):i\in\ZZ^d)$ of the type discussed here.
Examples appear 
in \citep{reza02, sepp00bd, sepp-3dhamm}. 
No explicit invariant distributions
are known for multidimensional height models. 
The variational scheme  proves 
 that scaled  height processes converge to solutions
of Hamilton-Jacobi equations 
 as in Theorem \ref{hydrothm}.   But again one can only 
assert the existence of $f$ and $g$ instead of giving them
explicitly.  

Another direction of generalization is to let the weights
 $\{Y_{i,j}\}$ have distributions other than  exponential or geometric. 
The height process $\hgt(t)$  ceases to be Markovian 
but the last-passage model of
Figs. \ref{fig1} and \ref{fig2} makes sense. 
As mentioned, the limit $\gamma(x,y)$ in \eqref{Glim} 
is explicitly known only for the  exponential and geometric cases.
  A distribution as simple as
Bernoulli ($Y_{i,j}$ takes only two values) cannot be handled. 
However, if the paths are altered to require that one or both
coordinates increase strictly, then the variational approach
 does find the explicit shape  for the  Bernoulli case 
\citep{sepp97incr, sepp98perc}.  

Thus the present situation is that an explicit limit shape
can be found 
only for some fortuitous combinations of path geometries
and weight distributions. 

\medskip

{\sl (ii) Partially asymmetric models.} 
Let us next address the case where the column variables 
$\hgti_i$ are allowed to jump both up and down.  Fix two
parameters $0<q<p$ such that $p+q=1$ (convenient normalization).
Give each column $i$ two independent Poisson clocks,
$N_i^{(-)}$ with rate $p$ and $N_i^{(+)}$ with rate $q$. 
At jump times of $N_i^{(\pm)}$ variable $\hgti_i$ attempts to
jump to $\hgti_i\pm 1$, and as before, a jump is completed 
if its execution does not take the state out of the state space
$\stsp_1$.   For the increment
process this means that a 
 {\tt 10} pair becomes a {\tt 01} pair at rate $p$, and the opposite
move happens at rate $q$.  Bernoulli distributions \eqref{defnurho}
are still invariant.  
This increment process is the {\sl asymmetric simple exclusion
process} (ASEP).
In the same vein one can allow $K$ particles  per site
and talk about asymmetric $K$-exclusion. 

The envelope property of Lemma \ref{varlm} is now lost. 
An alternative approach 
from \citep{reza-01} utilizes compactness of
the random semigroups of the height process.
Limit points are
characterized as Hamilton-Jacobi semigroups
 via the Lions-Nisio theorem
\citep{lion-nisi}. Thereby Theorem \ref{hydrothm} is derived 
 for one-dimensional asymmetric $K$-exclusion.
For $K=1$ the flux $f$ in \eqref{hj1} must be replaced by 
$f(\rho)=(p-q)\rho(1-\rho)$, while for $2\le K<\infty$ the flux is
unknown.
In the multidimensional case 
 it is not known  if the resulting equation
itself is random or not.

\section{Fluctuations}
\label{fluctsec}

A simple  way to create initial height functions $\hgt^n(0)$ that
satisfy assumption \eqref{llninit} is to take independent 
increments with distributions 
\[\bP[\eta^n_i(0)=1]=n\bigl(u_0(\tfrac{i}n)-u_0(\tfrac{i-1}n)\bigr),  \] and 
at the origin  assign the deterministic value 
$\hgti^n_0(0)=[nu_0(0)]$. The stationary situation is of 
this type with $u_0(x)=\rho x$.   Then initial fluctuations 
 \[ n^{-1/2}\{\hgti^n_{[nx]}(0)-nu_0(x)\} \]
are Gaussian in the limit $n\to\infty$.  This makes it 
 natural to look for  a  distributional limit at later times $t>0$ 
 on the central limit scale $n^{1/2}$:  
\beq
{n^{-1/2}}\{\hgti^n_{[nx]}(nt)- nu(x,t)\}
\longrightarrow \zeta(t,x) \quad\text{as $n\to\infty$,} 
\label{fluct2}\eeq
 for some limit process $\zeta(t,x)$.   Such limits can be proved
but process $\{\zeta(t,x)\}$ turns out to be a deterministic function 
of the initial fluctuations $\{\zeta(0,x)\}$. Consequently limit 
\eqref{fluct2} does not record any fluctuations created by the 
dynamics.  Theorem \ref{fluctthm2} below gives a precise statement
of this type.

In  asymmetric systems 
the fluctuations created by the dynamics occur on a scale smaller
than $n^{1/2}$.  Two types of such phenomena have been found.  
Processes related to the last-passage model and exclusion process
discussed in Section \ref{llnsec} have order $n^{1/3}$ fluctuations 
whose limits are distributions from random matrix theory. 
A class of linear processes has order $n^{1/4}$ fluctuations and 
Gaussian limits related to fractional Brownian motion with
Hurst parameter $H=1/4$.  
To see these lower order fluctuations one can start the system 
with a deterministic initial state, or one can start the 
system in the stationary distribution or some other random 
state but then  follow the evolution
along characteristic curves of the macroscopic PDE.  
The fluctuation  situation is very different for 
symmetric systems; the reader can consult
 \citep[Chapter 11]{kipn-land}.

\subsection{Exclusion process}
\label{asepsec}
Probability distributions from random matrix theory were
discovered as limit laws for last-passage growth models almost
a decade ago.

\begin{theorem} {\rm\citep{joha}.}  For the corner
growth model 
\beq
\bP\Bigl[\; \frac{G([xn], [ny]) -n\gamma(x,y)}{c(x,y)n^{1/3}}
\leq s\;\Bigr] \longrightarrow F(s)
\quad\text{as $n\to\infty$,}\label{joha1}\eeq
where   $F$ is the   Tracy-Widom GUE distribution. 
\label{johathm}\end{theorem}

The distribution $F$   first appeared as 
the limit distribution of the scaled largest eigenvalue of a random
Hermitian matrix from the GUE \citep{trac-wido-94}. 
GUE is short for {\sl Gaussian Unitary Ensemble}. This means that
a random Hermitian matrix is constructed by putting IID
complex-valued Gaussian random variables above the diagonal,
  IID
real-valued Gaussian random variables on the diagonal, and 
letting the Hermitian property determine the entries below
the diagonal.
Then as the matrix grows in size, the variances of the 
entries are scaled appropriately to obtain limits. The standard reference
is \citep{meht}.

Theorem \ref{johathm} and related results initially
arose  entirely outside probability theory (except for the
statements themselves), involving the RSK correspondence
and Gessel's identity from combinatorics and techniques
from integrable systems to analyze the asymptotics of the 
resulting determinants. 
The {\sl RSK correspondence}, named after Robinson, Schensted and 
Knuth, is a bijective mapping between certain arrays of
integers or integer matrices (in this case  the matrix in
Figure \ref{fig1} if the $Y_{i,j}$'s are integers) 
 and pairs of Young tableaux.  These latter objects are
ubiquitous in combinatorics. Standard references are 
\citep{fult, saga}.  
 More recently determinantal 
point processes have appeared as the link between the
growth processes and random matrix theory. 
We shall not pursue these topics further for 
many excellent reviews are available: \citet{baik-05},
\citet{deif-00}, \citet{joha-02}, \citet{koni-05} and  \citet{spoh-06}.

Precise limits such as 
\eqref{joha1}  have so far been restricted to
totally asymmetric systems.  Next  we discuss ideas that
fall short of exact limits but do give the correct order
of the variance of the height for partially asymmetric systems. 

Consider  the height process $\hgt(t)$  whose increments
$\eta_i(t)=\hgti_i(t)-\hgti_{i-1}(t)$ 
 form  the asymmetric simple exclusion process (ASEP). 
This process was introduced in the remarks at the end of 
Section \ref{llnsec}.  Each height variable  $\hgti_i$ attempts
 downward jumps with rate $p$  and upward jumps with rate $q$,
and $p>q$.  
A jump is suppressed if it would 
lead to a violation of the restrictions 
$\hgti_{i-1}\le \hgti_i$ or $\hgti_i\ge \hgti_{i+1}-1$ 
encoded in the state space $\stsp_1$ of \eqref{tasepX}. 
In the increment process each 
 {\tt 10} pair becomes a {\tt 01} pair at rate $p$ and each
  {\tt 01} pair becomes {\tt 10} pair at rate $q$.

On large space and time scales 
the height process obeys the Hamilton-Jacobi equation 
\[
\text{$u_t +  f(u_x) = 0$
with $f(\rho)=(p-q)\rho(1-\rho)$.}
\]  This PDE carries information
along the  curves  $\dot{x}= f'(u_x(t,x))$, in the sense that the 
slope $u_x$ is constant along these curves as long as it is 
continuous. 
At constant slope $u_x=\rho$ 
 the characteristic  speed is 
$ V^\rho=f'(\rho)=(p-q)(1-2\rho).$

Consider  the stationary process:
$0<\rho<1$ is fixed, and  at each time $t\ge 0$ 
the increments $\{\eta_i(t)\}_{i\in\ZZ}$ have 
Bernoulli $\nu_\rho$-distribution from \eqref{defnurho}. 
  Normalize the heights by setting initially $\hgti_0(0)=0$.  
We determine the order of magnitude of 
 the variance of the height as seen by an observer
traveling at speed $V^\rho$.  

\begin{theorem}
{\rm\citep{bala-sepp2}}  Height fluctuations along 
the characteristic satisfy 
\[
0< \liminf_{t\to\infty} {t^{-2/3}}{\bvar\{ \hgti_{[V^\rho t]}(t)\}}
\leq \limsup_{t\to\infty} {t^{-2/3}}{\bvar\{ \hgti_{[V^\rho t]}(t)\}}
<\infty.
\] 
\label{13thm} \end{theorem}

If the observer choses any other speed  $v\ne V^\rho$, only 
a translation of initial Gaussian fluctuations would be observed.
Take $v>V^\rho$ to be specific. 
Due to the normalization $\hgti_0(0)=0$
 we can write 
\beq
 \hgti_{[vt]}(t)=\bigl(\hgti_{[vt]}(t)-\hgti_{[(v-V^\rho)t]}(0)\bigr)
+\sum_{i=1}^{[(v-V^\rho)t]}\eta_i(0).
\label{fluct3}\eeq 
On the right the first expression in parentheses is a height
increment along a characteristic and so by the theorem has fluctuations 
of order $t^{1/3}$.  
The last sum of initial increments has  Gaussian fluctuations of order
 $t^{1/2}$ and consequently drowns out the first term. 

The proof of Theorem \ref{13thm} is entirely different from the 
proofs of Theorem \ref{johathm}.   As it involves an important 
probabilistic idea let us discuss it briefly.  

Couplings enable us to study the evolution of discrepancies
between processes.  In exclusion processes these discrepancies
are  called  
 second class particles.  Consider two initial ASEP configurations 
$\eta(0), \zeta(0)\in\{0,1\}^\ZZ$.  The configurations differ at the 
origin: 
 $\zeta(0)$ has a particle at the origin ($\zeta_0(0)=1$) but 
 $\eta(0)$ does not ($\eta_0(0)=0$). 
At all other sites $i\ne 0$ we give the configurations
a common but random value  
$\eta_i(0)=\zeta_i(0)$ according to the mean $\rho$ Bernoulli distribution. 
Let the joint process $(\eta(t),\zeta(t):t\ge 0)$ evolve
together  governed by the {\sl same}  Poisson  clocks. 
The effect of this coupling is that there is always exactly
one site $Q(t)$ such that $\zeta_{Q(t)}(t)=1$, $\eta_{Q(t)}(t)=0$, 
and $\eta_i(t)=\zeta_i(t)$ for all $i\ne Q(t)$.

$Q(t)$ is the location of a {\sl second class particle} relative
to the process $\eta(t)$.  (Relative to $\zeta(t)$ one should say
``second class {\sl anti}particle.'')  In addition to ordinary
exclusion jumps,  
 $Q$ yields to $\eta$-particles:  if an $\eta$-particle
at $Q+1$ jumps left (rate $q$) then $Q$ jumps right to switch places
with the $\eta$-particle. Similarly an $\eta$-particle at $Q-1$
switches places with $Q$ at rate $p$.  These special 
jumps follow from considering the effects of clocks $N^{(\mp)}_{Q\pm 1}$
on  the discrepancy between $\eta$ and $\zeta$. 

Proof of Theorem \ref{13thm} utilizes couplings 
of several  processes with different initial
conditions.  Evolution of 
 second class particles is directly related to differences in 
particle current (height)
between processes.  On the other hand $Q$ and the height variance
are related through this identity:
\beq
\bvar\{ h_{[vt]}(t)\} = \rho(1-\rho) 
\bE \bigl\{ \,\lvert Q(t)-[vt]\,\rvert\,\bigr\}
\quad\text{for any $v$.} 
\label{Qvarh}\eeq
The right-hand side can be expected to have order smaller than $t$ 
precisely when $v=V^\rho$ on account of this second identity: 
\beq
\bE Q(t)=tV^\rho.  \label{EQ(t)}\eeq
  From these ingredients  the bounds
in Theorem \ref{13thm} arise.

\medskip

{\bf Further remarks.} 
As already suggested at the end of Section \ref{llnsec},
a major problem for  growth models  is to find robust techniques
 that are not dependent on particular choices
of probability distributions or path geometries. 
 Progress on fluctuations of the corner growth model beyond the 
exponential case has come in situations that are
in some sense extreme: for distributions with heavy 
tails \citep{hamb-mart-07}
or for points close to the boundary of the quadrant
 \citep{baik-suid, bodi-mart}.  See
review by \citet{mart-06}.

The second class particle appears in many places in interacting
particle systems.  In the hydrodynamic limit picture
second class particles converge to characteristics and shocks of the 
macroscopic PDE \citep{ferr-font-94b, reza-95, sepp01IIclass}. 
Versions of identities \eqref{Qvarh} and \eqref{EQ(t)}
 are valid for zero range and bricklayer processes
 \citep{bala-sepp-07JSP}.  Equation \eqref{EQ(t)} is
 surprising because the process as seen by the second class
particle is {\sl not}  stationary. 

In general the view of the process from the second class particle 
is complicated. 
Studies  of  invariant distributions seen by  second
class particles appear in 
\citep{derr-etal-93, ferr-font-koha, ferr-mart-07}.   
There are special cases of parameter
values for certain processes where unexpected simplification 
takes place  and the process seen from the second class particle 
 has a product-form invariant distribution \citep{derr-lebo-spee-97,
bala-01}.

\subsection{Hammersley process}
\label{hammsec}
We began this paper with the exclusion process because this process
is by far the most studied among its kind. 
It behooves us to introduce also Hammersley's
process for which several important results were proved first.
  It has an elegant graphical construction 
that is related to a classical combinatorial question, namely
the maximal length of an increasing subsequence of a random
permutation.  This goes back to \citep{hamm}; see also 
\citep{aldo-diac95, aldo-diac99}.

\begin{figure}[ht]
\begin{center}
\begin{picture}(280,140)(5,20)
\put(40,30){\vector(1,0){200}}
\put(40,30){\vector(0,1){130}}
\put(80,60){\line(1,0){120}} \put(80,60){\line(0,1){80}}
\put(80,140){\line(1,0){120}}\put(200,140){\line(0,-1){80}}
\put(80,20){$a$} \put(200,20){$b$}
\put(30,58){$s$}\put(30,138){$t$}
\dottedline{3}(40,60)(80,60) \dottedline{3}(40,140)(80,140)
\dottedline{3}(80,30)(80,60) \dottedline{3}(200,30)(200,60)
\put(90,67){$\times$} \put(120,77){$\times$}
\put(130,100){$\times$} \put(180,125){$\times$}
\put(94.08,69.5){\circle{9}}\put(124.08,79.5){\circle{9}}
\put(134.08,102.5){\circle{9}} \put(184.08,127.5){\circle{9}}
\put(160,90){$\times$}\put(100,120){$\times$}\put(110,110){$\times$}
\put(150,40){$\times$}\put(162,47){$\times$}\put(60,49){$\times$} 
\put(210,137){$\times$}\put(220,100){$\times$}
\end{picture}
\end{center}
\caption{Increasing sequences among planar Poisson
points marked by $\times$'s. 
$L((a,s),(b,t))=4$ as shown by the circled Poisson points that 
form an increasing sequence.}
\label{fig:ulam}
\end{figure}
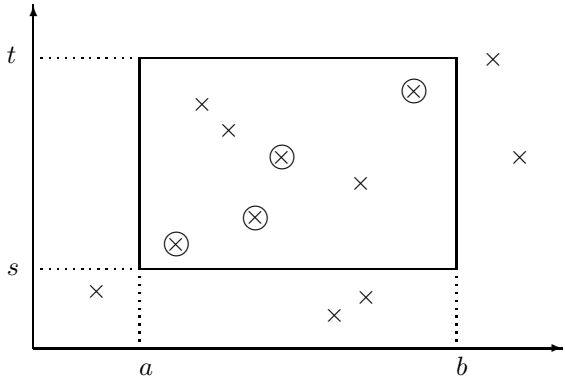

We begin with the growth model. 
Put a  homogeneous rate 1 Poisson point process on the plane.
This is a random discrete subset of the plane characterized by the 
following property:  the number of points in a Borel set $B$ 
is Poisson distributed with mean given by the area of $B$ and
independent of the points outside $B$.    
Call a sequence $(x_1,t_1), (x_2,t_2), \dotsc,(x_k,t_k)$  of these
Poisson  points 
{\sl increasing} if $x_1<x_2<\dotsm<x_k$   and  $t_1<t_2<\dotsm<t_k$.
Let $L((a,s),(b,t))$ be the maximal number of points on an 
increasing sequence in the rectangle $(a,b\,]\times(s,t]$
(Fig.\ \ref{fig:ulam}).  
The random permutation comes from mapping the 
ordered $x$-coordinates to ordered 
$t$-coordinates in the rectangle, and 
$L((a,s),(b,t))$ is precisely the maximal length of an 
increasing subsequence of this permutation.

The limit  
\beq \lim_{n\to\infty} n^{-1}L((0,0),(nx,nt))=
2\sqrt{xt} \quad\text{a.s.} \label{ulamlim}\eeq
holds for $x,t>0$.  The limit exists by superadditivity 
exactly as for \eqref{Glim}.  The functional form 
$c\sqrt{xt}$ follows from scaling properties of the 
Poisson process.  
 The value  $c=2$ was first derived by \citet{vers-kero-77},
while \citet{loga-shep-77} independently proved $c\ge 2$. 
The fluctuation result  for $L$ is analogous to \eqref{joha1},
with normalization $n^{1/3}$ and  the Tracy-Widom
limit \citep{baik-deif-joha-99}.

\begin{figure}[ht]
\begin{center}
\begin{picture}(280,140)(5,10)
\put(50,20){\line(1,0){230}}
\put(50,120){\line(1,0){230}}
\put(12,17){time $0$} \put(12,117){time $t$}
\put(97.5,17){\large$\bullet$} \put(86,9){$z_{i-1}(0)$}
\put(100,20){\line(0,1){30}} \put(100,50){\line(-1,0){15}}
\put(85,50){\line(0,1){70}} \put(81.3,48){$\times$}
\put(81.5,117){\large$\bullet$} \put(70,127){$z_{i-1}(t)$}
\put(157.5,17){\large$\bullet$} \put(150,9){$z_{i}(0)$}
\put(160,20){\line(0,1){50}} \put(160,70){\line(-1,0){30}}
\put(130,70){\line(0,1){10}} \put(126.5,68){$\times$}
\put(130,80){\line(-1,0){10}} 
\put(120,80){\line(0,1){40}} \put(116.3,77.8){$\times$}
\put(116.5,117){\large$\bullet$} \put(109,127){$z_{i}(t)$}
\put(247.5,17){\large$\bullet$} \put(235,9){$z_{i+1}(0)$}
\put(250,20){\line(0,1){15}} \put(250,35){\line(-1,0){70}}
\put(180,35){\line(0,1){85}} \put(176.5,32.7){$\times$}
\put(176.5,117){\large$\bullet$} \put(163,127){$z_{i+1}(t)$}
\put(220,85){\line(0,1){35}\line(1,0){60}} \put(216.3,82.8){$\times$}
\put(216.5,117){\large$\bullet$} \put(207,127){$z_{i+2}(t)$}
\end{picture}
\end{center}
\caption{Portion of the graphical construction of Hammersley's process.
$\times$'s mark space-time Poisson points. $\bullet$'s mark particle
locations at time $0$ and at a later  time $t>0$.  Space-time 
trajectories of particles 
are shown. The horizontal segments are traversed instantaneously
and the vertical segments at constant speed 1.}  \label{figHamm}
\end{figure}
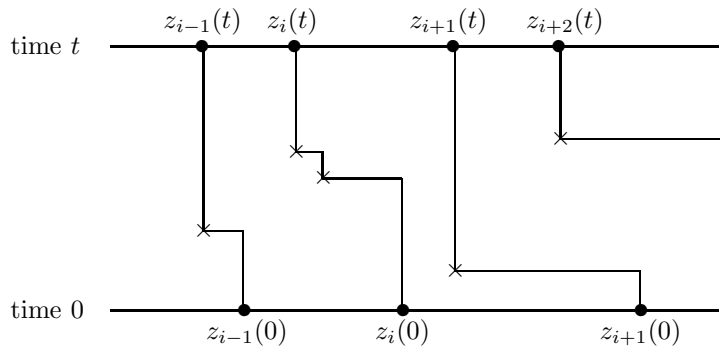

We embed the increasing sequences in the graphical
construction of the {\sl Hammersley
process}.   This process consists of point particles
that move on  $\RR$ by jumping. 
  Put a rate 1 Poisson point process
on the space-time plane    and place the particles initially 
on the real axis.  Move the real axis up
 at constant speed 1. 
Each Poisson point $(x,t)$ instantaneously  pulls to $x$ 
 the next particle to the right of $x$. 
We label the particles from left to right:
$\hham_i(t)\in\RR$ is the position of particle $i$ at time $t$. 
We could regard the variables $\hham_i$ as heights again,
but the particle picture seems more compelling. 
This construction is illustrated by Fig.\ \ref{figHamm}.
In terms of infinitesimal rates, the construction realizes this rule:
independently of other particles, at rate $z_i-z_{i-1}$
variable $z_i$ jumps to a uniformly chosen location in the
interval $(z_{i-1},z_i)$.

  As in Section \ref{llnsec}, there
is a variational characterization for  this construction. 
Define an inverse for the maximal path variable 
$L((a,s),(b,t))$  by
\[
\Gamma((a,s),t,w)=\inf\{h\ge 0: L((a,s),(a+h,t))\ge w\}.
\] 
Take an initial particle 
  configuration  $\{\hham_i(0)\}\in\RR^\ZZ$ that satisfies 
$\hham_{i-1}(0)\le \hham_i(0)$ and 
$i^{-2}\hham_i(0)\to 0$ as $i\to-\infty$.  
Then the graphical construction leads to a well-defined 
evolution $\{\hham_i(t)\}$ that satisfies
\beq
\hham_i(t)=\inf_{k: k\le i}\bigl\{ \hham_k(0)+
\Gamma((\hham_k(0),0),t,i-k)\bigr\}.
\label{varhamm}\eeq 

Here is the hydrodynamic limit.
Consider  a sequence of processes $\hham^n(t)$ 
 that satisfies 
$n^{-1}\hham^n_{[ny]}(0)\to u_0(y)$ for each $y\in\RR$,
say in probability.  The initial function
 $u_0$ is  nondecreasing, locally Lipschitz and satisfies
\beq y^{-2}u_0(y)\to 0\quad\text{ as $y\to-\infty$.} \label{u0ass}\eeq 
  Define 
\beq
u(t,x)=\inf_{y: y\le x}\Bigl\{ u_0(y)+\frac{(x-y)^2}{4t}\Bigr\},
\quad (t,x)\in(0,\infty)\times\RR.
\label{hopflaxhamm}\eeq
Since rates are unbounded, we need to assume a left tail bound
to prevent the particles from disappearing to $-\infty$:
given $\e>0$ there exist $0<q, n_0<\infty$ such that   
\[
\bP\{\text{$\hham^n_i(0)<-\e i^2/n$ for some $i\le -qn$}\} \le \e
\quad\text{for $n\ge n_0$.} 
\]
 Under these  assumptions 
\beq
n^{-1}\hham^n_{[nx]}(nt)\to u(t,x) \quad\text{ in probability }
\label{hammlln}\eeq
 \citep{sepp-96}.
The function defined by the Hopf-Lax formula 
\eqref{hopflaxhamm}  solves the Hamilton-Jacobi
equation $u_t+(u_x)^2=0$.  

Let us state a precise 
result about the central limit scale fluctuations
\eqref{fluct2}  that covers also shocks.  For $(t,x)\in(0,\infty)\times\RR$ 
  let 
\beq
I(t,x)=\Bigl\{y\in(-\infty, x]: 
u(t,x)= u_0(y)+\frac{(x-y)^2}{4t}\Bigr\}
\label{defIxt}\eeq
be the set of minimizers in \eqref{hopflaxhamm}, guaranteed
nonempty and compact by hypothesis \eqref{u0ass}.
Then $(t,x)$ is a {\sl shock} if $I(t,x)$ is not a singleton. This is
equivalent to the nonexistence of the $x$-derivative $u_x(t,x)$. 

Fluctuations on the scale $n^{1/2}$ from the limit 
\eqref{hammlln}  are described by the process
\[
\zeta_n(t,x)= {n^{-1/2}} \{\hham^n_{[nx]}(nt)- nu(x,t)\}.
\]
Assume the existence of a continuous random function $\zeta_0$ on $\RR$ 
such that the convergence in distribution 
$
\zeta_n(0,\,\cdot\,)\to \zeta_0
$
holds in the topology of uniform convergence on compact sets.
Define the process $\zeta$ by
\[
\zeta(t,x)=\inf_{y\in I(t,x)}\zeta_0(y)
\]
where $I(t,x)$ is the (deterministic) set defined in \eqref{defIxt}. 

\begin{theorem}
For each $(t,x)$,  $\zeta_n(t,x)\to \zeta(t,x)$
in distribution.  
\label{fluctthm2}\end{theorem} 

As stated in the beginning of  Section \ref{fluctsec}, 
this distributional limit reflects no
contribution from dynamical fluctuations as the process
$\zeta$ is a deterministic transformation of $\zeta_0$.  
The underlying reason is that the dynamical fluctuations of order
$n^{1/3}$ are not visible on the $n^{1/2}$ scale.  
The dynamical fluctuations are the universal ones described 
by the Tracy-Widom laws. See again the discussion and 
references that follow Theorem \ref{johathm}.

\medskip

{\bf Further remarks.}  
The {\sl polynuclear growth model} (PNG) is another related
(1+1)-dimensional growth model used  by several authors 
 for studies of  Tracy-Widom fluctuations and the Airy process
in the 
KPZ scaling picture \citep{baik-rain-00, ferrPL-04,
joha-03, prah-spoh-02,
prah-spoh-04}. Like the Hammersley process,
the graphical construction of the PNG 
 utilizes a planar Poisson process, and in fact the same 
underlying last-passage model of increasing paths. 
This time the Poisson points mark space-time nucleation events from which
new layers grow laterally at a fixed speed.  Roughly speaking, this
corresponds to putting the time axis at a 45 degree angle
in Figs.\ \ref{fig:ulam} and \ref{figHamm}. 

More about the phenomena related to Theorem \ref{fluctthm2} 
 can be found in article \citep{sepp02diff}.  A similar  theorem for
TASEP  appears in
\citep{reza02tasep}.  Earlier work on the diffusive fluctuations
of ASEP was done by \citet{ ferr-font94a, ferr-font-94b}.

\subsection{Linear models}

We turn to systems macroscopically governed by  linear 
first order equations $u_t+bu_x=0$.  Fluctuations across the 
characteristic  occur now on the scale $n^{1/4}$ and converge to 
a Gaussian process 
related to fractional Brownian motion. 

\medskip

The {\sl random average process} (RAP) was first studied by 
\citet{ferr-font-RAP}.  The state of the process 
is a    height  function $\hrap:\ZZ\to\RR$ with $\hrap_i\in\RR$
denoting the height over site $i$.   (More generally
the domain can be $\ZZ^d$.) 
 The basic step
of the evolution is that a value $\hrap_i$ is replaced by a 
weighted average of values in a neighborhood, and the randomness comes
 in the weights. 
This time we consider a discrete time process. The basic step
is carried out simultaneously at all sites $i$.

Now  precise formulations. Let 
$\{u(k,\tau): k\in\ZZ, \,\tau\in\NN\}$ be an IID collection of random
probability vectors indexed by space-time $\ZZ\times\NN$.  
In terms of coordinates 
 $u(k,\tau)=(u_j(k,\tau): -M\leq j\leq M)$.  We assume the
system has finite range defined by the fixed parameter $M$. 
We impose a minimal assumption that guarantees that the
weight vectors are not entirely degenerate: 
\beq \bP\{ \max_j u_j(0,0) < 1\}>0. \label{rapass1}\eeq
For technical convenience we also assume that the process is
``on the correct lattice'': there does not exist
an integer $h\ge 2$ such that for some $b\in\ZZ$ the mean weights
 $p(j)=\bE u_j(0,0)$ satisfy $\sum_{j\in b+h\ZZ}p(j)=1$. 

To start the dynamics 
let $\hrap(0)$ be a given random or deterministic initial height function.
 The process $\hrap(\tau)$, 
 $\tau=0,1,2,\dotsc$,  is defined iteratively by 
\beq \hrap_i(\tau)=\sum_{j }u_j(i,\tau) \hrap_{i+j}(\tau-1)\,,
\quad \tau\ge 1, \,i\in\ZZ. \label{rapstep}\eeq 
As before, we can define 
the process of increments 
$\eta_i(\tau)=\hrap_i(\tau)-\hrap_{i-1}(\tau)$.  The increments
 also  evolve via random linear mappings and are conserved
like particles in exclusion processes.  

As in Section \ref{llnsec} we create suitable initial conditions
for a hydrodynamic limit.  Consider a  sequence of 
processes $\hrap^n(\tau)$  
indexed by $n\in\NN$, initially  normalized by  $\hrap^n_0(0)=0$.
For each $n$ assume 
independent  initial increments  
$\{\eta^n_i(0):i\in\ZZ\}$ with 
 \beq
\bE [\eta^n_i(0)]=\rho(i/n) 
\quad\text{and}\quad 
 \bvar [\eta^n_i(0)]=v(i/n)
\label{rapass2}\eeq
for given H\"older $1/2+\e$ functions  $\rho$ and $v$. 
Assume a uniform   moment bound:
$\sup_{n.i}\bE[\,\lvert \eta^n_i(0)\rvert^{2+\delta}\,]<\infty$ 
for some $\delta>0$.

 The hydrodynamic limit  is  rather 
trivial for it consists only of translation. 
Define a function $u$ on $\RR$ by $u(0)=0$ and $u'(x)=\rho(x)$. 
The characteristic speed is 
\[ b=-\sum_j j p(j).\]
Then for each $(t,x)\in\RR_+\times\RR$, 
\[ n^{-1}\hrap^n_{[nx]}([nt]) \longrightarrow u(x-bt)
\qquad \text{as $n\to\infty$, in probability.} \]
In other words,  the  height obeys the  linear PDE
$u_t+ bu_x=0.$

On the central limit scale $n^{1/2}$ one would also see only
 translation of initial fluctuations. 
To see something nontrivial we look at 
 fluctuations around a characteristic line. 
Fix a point  $\ybar\in\RR$ and consider the 
 characteristic line $t\mapsto \ybar+tb$ emanating  from $(\ybar,0)$.
Define space-time process
 \[Z_n(t,r)= 
 \hrap^n_{[n\ybar]+[r\sqrt{n}\,]+[ntb]}([nt])\,-\,
 \hrap^n_{[n\ybar]+[r\sqrt{n}\,]}(0)\,,\quad
 (t,r)\in\RR_+\times\RR. 
\]
The  spatial variable $r$
 describes fluctuations around the characteristic on the
spatial scale $n^{1/2}$.  
For the  increment process  $Z_n(t,0)$ represents 
 net current  from right to left across the
characteristic.

\begin{theorem}
{\rm \citep{bala-rass-sepp}} 
The finite-dimensional distributions of the process 
$n^{-1/4}Z_n $ converge to those of the Gaussian process
 $\{z(t,r): t\geq 0, \,r\in\RR\}$ described below. 
\label{rapthm2}\end{theorem} 

The statement means that for any finite collection of space-time points
$(t_1,r_1),\dotsc,(t_k,r_k)$, the $\RR^k$-valued random vector
$n^{-1/4}(Z_n(t_1,r_1),\dotsc,Z_n(t_k,r_k))$ converges in 
distribution to the
vector $(z(t_1,r_1),\dotsc, z(t_k,r_k))$.  
The limiting process $z$ has the following representation in terms
of stochastic integrals: 
\beq \begin{split}
z(t,r)&\ =\ \rho(\ybar)\sigma_a\sqrt{\kappa}\iint_{[0,t]\times\RR} 
\varphi_{\sigma_a^2(t-s)}(r-z)\,dW(s,z)\\
&\qquad +\   \sqrt{v(\ybar)}\;\int_\RR 
\sign(x-r)\,
\Phi_{\sigma_a^2t}\bigl(-\lvert{x-r}\rvert\,\bigr)\,dB(x).
\end{split}
\label{defz}\eeq
Above  $W$ is a 2-parameter Brownian motion on 
$\RR_+\times \RR$ and 
  $B$ is a  1-parameter Brownian motion 
on $\RR$ independent of  $W$.  The first integral represents
 dynamical noise generated by the random weights, and the second
the initial noise propagated by the evolution.  The functions
in the integrals are Gaussian densities and distribution functions:
\[
\varphi_{\sigma^2}(x)=\frac1{\sqrt{2\pi \sigma^2}}
\exp\Bigl\{-\frac{x^2}{2\sigma^2} \Bigr\}\quad\text{and}\quad  
\Phi_{\sigma^2}(x)=\int_{-\infty}^x \varphi_{\sigma^2}(y)\,dy.
\]
The only effects  from the  initial height 
are the mean
$\rho(\ybar)$  and
variance $v(\ybar)$ of the increments around the point $n\ybar$. 
The parameter $\sigma_a^2$ is the variance of the probabilities
$p(j)$, and $\kappa$ another parameter determined by the 
distribution of the weights. 
Process $z$ has a self-similarity property: 
$\{z(at, {a}^{1/2}r)\}\overset{d}=\{a^{1/4}z(t,r)\}$.

In the special case where $v(\ybar)=\kappa\rho(\ybar)^2$
the temporal process 
 $\{z(t,r):t\in\RR_+\}$ (for any fixed $r$) 
has covariance 
\[
\bE z(s,r)z(t,r)=\frac{\sigma_a\kappa\rho^2}{\sqrt{2\pi}}
\bigl( \sqrt{s}+\sqrt{t}-\sqrt{\lvert{t-s}}\rvert\,\bigr).
\]
This identifies $z(\cdot\,,r)$
 as {\sl fractional Brownian motion} with Hurst parameter $H=1/4$. 
In particular, this limit arises  in a  stationary case
 where  the averaging involves two points and the weight is beta
distributed \citep[Example 2.1]{bala-rass-sepp}. 

The proof of Theorem \ref{rapthm2} utilizes a special case of another
stochastic model of great contemporary interest, namely 
random walk in random environment (RWRE). 
Here is how the RWRE arises.  An environment
$\omega=\{u(k,\tau)\}$  is determined by
the weight vectors.
Given $\omega$, define   a ``backward''  
walk   $\{X^{i,\tau}_s: s\in \ZZ_+\}$  
 on $\ZZ$ with initial position 
$X^{i,\tau}_0=i$ and transition probability 
\[ 
P^\omega(X_{s+1}^{i,\,\tau}=y\,|\,X^{i,\,\tau}_s=x)=u_{y-x}(x, \tau-s),
\quad s=0,1,2,\dotsc
\]
The superscript $\omega$ on  
$P^\omega$ indicates that it is the {\sl quenched}  path measure 
of $\{X^{i,\tau}_s:s\in\ZZ_+\}$ under a fixed $\omega$. 
The basic step \eqref{rapstep}  of 
RAP evolution can be rewritten so that $\hrap_i(\tau)$ 
 equals  the average value of the previous 
height function $\hrap(\tau-1)$ seen by a walk started at $i$
 after one step: 
\[ 
\hrap_i(\tau)=\sum_j 
u_{j-i}(i,\tau)  \hrap_j(\tau-1)
=E^\omega\bigl[ \hrap_{X_{1}^{i,\,\tau}}(\tau-1)\bigr].
\] 
This can be iterated all the way down to the initial height function:
\[ 
\hrap_i(\tau)=E^\omega\bigl[\hrap_{X^{i,\,\tau}_\tau}(0)\bigr].
\]
 Note that the expectation  $E^\omega$ over 
paths of the walk $X^{i,\,\tau}_s$ under fixed weights $\omega$ 
 sees the  initial height function
  $\{\hrap_i(0)\}_{i\in\ZZ}$ as a    constant. 

We have here a special type of RWRE called ``space-time.'' 
Another term used is
 ``dynamical environment'' because after each step the walk 
sees a new sample of its environment.  Proof of Theorem \ref{rapthm2}
requires limits for the walk itself and its quenched
mean process $E^\omega(X^{i,\,\tau}_s)$. These results
appear in \citep{bala-rass-sepp, rass-sepp1}.  

\medskip

Independent walks on $\ZZ$ display  the same behavior as RAP. 
Let the process $Z_n(t,r)$ be the net particle
current across the characteristic 
$t\mapsto [n\ybar]+[r\sqrt{n}\,]$ $+$ $[ntb]$ where $b$ is
 the common average speed of the particles.  Then
under suitable assumptions on the initial particle 
arrangements and their jump kernel,  $Z_n$ satisfies
a stronger form of Theorem \ref{rapthm2} that also
contains 
process-level convergence.  One adjustment is necessary: the 
 constants in front of the stochastic integrals in \eqref{defz}
are different for the random walk case.  
The stationary system sees again fractional 
Brownian motion as the limit of the current. 
Details for the random walk case appear in 
\citep{sepp-rw, kuma-07}. 
Earlier related results for a Poisson system of independent Brownian
motions appeared in \citep{durr-gold-lebo}.

\section{Large deviations}
\label{ldpsec}

We present the large deviation picture for the 
Hammersley process, so we continue 
 in the setting of Section \ref{hammsec}. 
Recall the definition of the longest path model among planar 
Poisson points illustrated in Fig.\ \ref{fig:ulam}. 
Abbreviate $L_n=L((0,0),(n,n))$. Then the limit \eqref{ulamlim} is
  $n^{-1}L_n\to 2$. 
Here is the large deviation theorem for $L_n$.  It was 
completed shortly before the fluctuation result
of \citet{baik-deif-joha-99}, through
a combination of several independent papers: \citet{loga-shep-77},
\citet{kimj-96},  \citet{sepp-98ldp} and \citet{deus-zeit}.

\begin{theorem} 
We have the following upper and lower tail large deviation bounds.
\beq
\lim_{n\to\infty} n^{-1}\log \bP\{L_n\ge nx\} =-I(x)
\quad\text{for $x\ge 2$ }
\label{Lutldp}\eeq
with rate function
$
I(x)=2x\cosh^{-1}(x/2)-2\sqrt{x^2-4\,}.
$
\beq
\lim_{n\to\infty} n^{-2}\log \bP\{L_n\le nx\} =-U(x)
\quad\text{for $0\le x\le 2$}
\label{Lltldp}\eeq
with rate function
$U(x) =\int_x^2 R_2(s)\,ds$ 
where $R_2(s)=s\log(s/2)-s+2 $ is the rate function for IID
mean 2 Poisson random variables.
\label{ldpthm}\end{theorem} 


To develop this theme further
 we state a lower tail LDP for the tagged particle in Hammersley's
process.  An interesting feature is  that the large deviation
 rate functions
again obey  the Hopf-Lax semigroup formula, as did the limit
\eqref{hammlln}. 
 The assumption is that lower tail rate functions exist 
initially: for all $y,s\in\RR$ the limit 
\[
J_0(y,s)=-\lim_{n\to\infty} n^{-1}\log \bP\{\hham^n_{[ny]}(0)\le ns\}
\]
exists and is left continuous in $y$ for each fixed $s$.  Define 
\beq
\Psi(w,r)=-\lim_{n\to\infty} n^{-1}\log 
\bP\{\Gamma((0,0),(n,nw))\le nr\}.
\label{defPsi}\eeq
This limit exists by superadditivity. 
For technical reasons  a uniform  tail
bound is needed for the  initial particle locations: there exist
constants $0<C_j<\infty$ such that 
\[
\bP\{ \hham^n_i(0)\le -C_1\abs{i}\,\} \le e^{-C_2\abs{i}}
\quad\text{ for $i\le -C_3n$, for large enough $n$.} \] 

\begin{theorem} {\rm \citep{sepp-98ldp}}  The limit 
\[
J_{t}(x,r)=-\lim_{n\to\infty} n^{-1}\log \bP\{\hham^n_{[nx]}(nt)\le nr\}
\]
exists for all $x,r\in\RR$ and $t>0$, and is given by 
\beq
J_t(x,r)=\inf_{(y,s):y\le x,\,s\le r}
\Bigl\{ J_0(y,s)+t\Psi\Bigl(\frac{x-y}t, \frac{r-s}t\Bigr)\Bigr\}.
\label{ldphopflax}\eeq 
\label{ldpthm2}\end{theorem}

The approach  of Section \ref{llnsec} can be adapted to 
prove the upper tail LDP  \eqref{Lutldp} and 
Theorem \ref{ldpthm2}.  The stationary systems make 
explicit calculation again possible. 
Presently it is not clear how to include the
 lower tail LDP \eqref{Lltldp} in the variational 
framework.  Hence this part  requires a separate proof.
Details appear in \citet{deus-zeit} and \cite{sepp-98ldp}.

\medskip

{\bf Further remarks.}
Even though Theorem \ref{ldpthm} has explicit rate functions, 
this large deviation problem remains unfinished in an important
sense.  It is not understood how
 the system behaves to create a deviation, and it is not
clear what the rate functions $I$ and $U$ represent.  
The present proofs are too indirect.   

Let us illustrate 
through the random walk LDP \eqref{rwldp1}--\eqref{rwldp2}
how a large deviation problem ideally should be understood. 
To create a deviation $S_n\approx nu$ with $u> v$, the 
entire walk behaves as a random walk with mean step $u$.
Namely, it can be proved that the conditioned measure
$\bP(\,\cdot\,\vert\, S_n\ge nu)$ converges on the path
space  to the distribution
$\bQ^{(u)}$  of a mean $u$ random walk.  The value  $I(u)$ of the rate function
in \eqref{rwldp2} 
is the entropy of this measure  $\bQ^{(u)}$ relative to the original $\bP$.

Results of the type presented in this section appear
for TASEP in \citep{sepp98ebp}. 
The asymptotic analysis of \cite{baik-deif-joha-99} and
\cite{joha} gives also LDP's for the growth models. 
  Concentration results 
for a Brownian last-passage model appear in 
\citep{hamb-mart-ocon} and a general discussion of deviation 
inequalities for growth models in the lectures of \cite{ledo-07}. 



\bibliographystyle{chicago}
\bibliography{growthrefs}
\end{document}